\documentclass[11pt]{article}
\usepackage{amsfonts}
\usepackage{mathrsfs}

\usepackage[latin1]{inputenc}
\usepackage{amsmath,amssymb}
\usepackage{latexsym}
\usepackage[active]{srcltx}
\usepackage{multicol}  
\usepackage{multirow} 
\usepackage{threeparttable}  
\usepackage[
bookmarks=true,         
bookmarksnumbered=true, 
colorlinks=true, pdfstartview=FitV, linkcolor=blue, citecolor=blue,
urlcolor=blue]{hyperref}

\newtheorem{theorem}{Theorem}[section]
\newtheorem{corollary}[theorem]{Corollary}

\numberwithin{equation}{section}
\def\qed{{\hfill $\square$ \bigskip}}
\def\R{{\mathbb R}}
\def\E{{{\mathbb E}\,}}

\def\N{{\mathbb N}}

\def\square{{\vcenter{\vbox{\hrule height.3pt
        \hbox{\vrule width.3pt height5pt \kern5pt
           \vrule width.3pt}
        \hrule height.3pt}}}}
\def \eref#1{\hbox{(\ref{#1})}}

\begin{document}

\title{Kernel entropy estimation for linear processes II}
\date{\today}
\author{Yudan Xiong and Fangjun Xu\thanks{F. Xu is partially supported by National Natural Science Foundation of China (Grant No.12371156).} \\
}
\maketitle

\begin{abstract}  

Let $X=\{X_n:  n\in \N\}$ be a linear process with  bounded probability density function $f(x)$. Under certain conditions, we use the kernel estimator
\[
\frac{2}{n(n-1)h_n} \sum_{1\le i<j\le n}K\Big(\frac{X_i-X_j}{h_n}\Big)
\]
to estimate the quadratic functional of $\int_{\R}f^2(x)dx$ of the linear process $X=\{X_n:  n\in \N\}$ and  improve the corresponding results in \cite{ssx}.
\vskip.2cm

\vskip.2cm \noindent {\bf Keywords}:
\noindent linear process, kernel entropy estimation, quadratic functional, projection operator


\vskip.2cm 

\end{abstract}

\section{Introduction}

Let $X=\{X_n:\, n\in\mathbb{N}\}$ be a linear process defined by 
\begin{equation}\label{def-Xn}
X_n=\sum_{i=0}^{\infty}a_{i}\varepsilon_{n-i},
\end{equation}
where the innovations $\varepsilon_i$ are independent and identical distributed (i.i.d.) real-valued random variables in some probability space $(\Omega, \mathcal{F},\mathbb{P})$ and $a_i$ are real coefficients such that $\sum\limits^n_{i=0} a_i\varepsilon_{n-i}$ converges in distributions. Assume that the linear process  $X=\{X_n:\, n\in\mathbb{N}\}$  has a bounded density function $f$. Clearly, the quadratic functional $\int_{\R} f^2(x)\,dx$ is well-defined and it plays an important role in the study of the  entropies of the linear process $X=\{X_n:\, n\in\mathbb{N}\}$, say, quadratic R\'{e}nyi entropy $R(f)=-\ln (\int_{\R} f^2(x)\,dx)$ and Shannon entropy $S(f)=-\int_{\R} f(x)\ln f(x)\,dx$, see \cite{ssx} and references therein. To estimate the quadratic functional $\int_{\R} f^2(x)\, dx$ of the linear process $X=\{X_n:\, n\in\mathbb{N}\}$ defined in $\eref{def-Xn}$, the authors in \cite{ssx} used the kernel estimator
\begin{align} \label{def-Tn}
T_n(h_n)=\frac{2}{n(n-1)h_n} \sum_{1\le j< i\le n}K\Big(\frac{X_i-X_j}{h_n}\Big),
\end{align}
where the kernel $K$ is a symmetric and bounded function with $\int_{\R} K(u)\, du = 1$ and $\int_{\R} u^2|K(u)|\, du<\infty$. The bandwidth sequence $h_n$ satisfies $0<h_n\to 0$ as $n\to\infty$. Then, under certain conditions, asymptotic properties of the estimator $T_n(h_n)$ were established, see Theorems 2.1 and 2.2 in \cite{ssx}.  In this paper, we are going to improve these two theorems. Some ideas are borrowed from \cite{lx}. 

Let $\phi_{\varepsilon}(\lambda)$ be the characteristic function innovations. That is, for each $\lambda\in\mathbb{R}$, $\phi_{\varepsilon}(\lambda)=\mathbb{E}[e^{\iota \lambda \varepsilon_{1}}]$ where $\iota$ is the imaginary unit $\sqrt{-1}$. Moreover, for any integrable function $g(x)$, its Fourier transform is defined as~$\widehat{g}(u)=\int_{\mathbb{R}}e^{\iota x u}g(x)\, dx$. The following is our main result.

\begin{theorem} \label{thm1} 
For some $\gamma\in(0,1]$, assume that 
\begin{enumerate}
\item[{\bf (A1)}] there are at least three coefficients $\{a_i: i\geq 0\}$ of the linear process $X$ are non-zero and $\sum\limits^{\infty}_{i=0}|a_i|^{\gamma}<\infty$, 
\item[{\bf (A2)}] there exists a positive constant $c_\gamma$ such that $
\E|e^{\iota \lambda \varepsilon_1}-\phi_{\varepsilon}(\lambda)|^{2}\leq c_{\gamma} \left(|\lambda|^{2\gamma}\wedge 1\right)$
for all $\lambda\in\R$,
\item[{\bf (A3)}] (i) $\int_{\R}|\lambda|^{3\gamma} |\phi_{\varepsilon}(\lambda)|^2\, d\lambda<\infty$  or (ii) $\int_{\R}|\lambda|^{2\gamma} |\phi_{\varepsilon}(\lambda)|^2\, d\lambda<\infty$ and $\sup\limits_{\lambda\in\R}|\widehat{K}(\lambda)||\lambda|^{\gamma}<\infty$.
\end{enumerate}
Then there exist positive constants $c_1$ and $c_2$ such that
\begin{align} \label{r4}
\Big|\E T_n(h_n)-\int_{\R} f^2(x)\, dx \Big|\leq c_1\Big(\frac{1}{n}+h^{2\gamma}_n\Big),
\end{align}
\begin{align}  \label{r5}
\E\Big|T_n(h_n)-\E T_n(h_n)-\frac{1}{n}\sum^n_{i=1}Y_i\Big| \leq c_2 \Big(\frac{1}{\sqrt{n^3h^2_n}}+\frac{1}{\sqrt{n^2h_n}}+\frac{h^{\gamma}_n}{\sqrt{n}}\Big),
\end{align}
where $Y_i=2\big(f(X_i)-\E f(X_i)\big)$ for $1\leq i\leq n$, and, if additionally $nh_n\to\infty$ as $n\to\infty$,
\begin{align} \label{r6}
\sqrt{n}\, \Big[T_n(h_n)-\E T_n(h_n)\Big]\overset{\mathcal{L}}{\longrightarrow} N(0,4\sigma^2)
\end{align}
for some $\sigma^2\in(0,\infty)$.
\end{theorem}

\noindent
{\bf Remark}: (i) If at most two coefficients of the linear process $X$ are non-zero, then $X$ is an i.i.d or $2$-dependent sequence. The corresponding results are available in \cite{gn,kls,ssx}.

(ii) By Assumptions {\bf (A1)} and {\bf (A3)}, we could easily obtain that $\mathbb{E}[e^{\iota u X_n}]=\prod\limits^{\infty}_{i=0}\phi_{\varepsilon}(u a_i)\in L^1(\mathbb{R})$. So the probability density function $f(x)$ of the underlying linear process $X$ is bounded.

Next we compare our result with those in \cite{ssx}. The assumption (3) 
\[
\E|e^{\iota \lambda \varepsilon_1}-\phi_{\varepsilon}(\lambda)|^{2}\leq c_{\gamma,4} \left(|\lambda|^{4\gamma}\wedge 1\right)
\]
in Theorem 2.1 of \cite{ssx} is weaken to our assumption {\bf (A2)}. According to Example 3.2 in \cite{ssx}, our Theorem \ref{thm1} works for more linear processes than Theorem 2.1 in \cite{ssx}. The result (9)
\begin{align*}   
\E\Big|T_n(h_n)-\E T_n(h_n)-\frac{1}{n}\sum^n_{i=1}Y_i\Big| \leq c_4 \Big(\frac{1}{nh_n}+\frac{h^{\gamma}_n}{\sqrt{n}}\Big)
\end{align*}
 in Theorem 2.2 of \cite{ssx} is improved to (\ref{r5}) in our Theorem \ref{thm1}. Moreover, to obtain the central limit theorem in (\ref{r6}), we only require the bandwidth to satisfy $\lim\limits_{n\to\infty}nh_n=\infty$ instead of $\lim\limits_{n\to\infty}\sqrt{n}h_n=\infty$ in Theorem 2.2 of \cite{ssx}.

As a consequence of Theorem \ref{thm1}, we can easily obtain the following result.

\begin{corollary}
Under the assumptions {\bf (A1)}, {\bf (A2)} and {\bf (A3)} in Theorem \ref{thm1}. Let $h_n$ have the order $
\left\{\begin{array}{ll}
n^{-\frac{3}{4\gamma+2}} &  \text{for}\quad \gamma\in(0,1/4],\\ \\
n^{-\frac{2}{4\gamma+1}} &  \text{for}\quad \gamma\in(1/4,1].
\end{array}\right.$ If $0<\gamma\leq 1/4$, then 
\[
T_n(h_n)-\int_{\mathbb{R}} f^2(x)\, dx=O_{\mathbb{P}}(n^{-\frac{3\gamma}{2\gamma+1}});
\]
if $1/4<\gamma\leq 1$, then 
\[
\sqrt{n}\Big[T_n(h_n)-\int_{\mathbb{R}} f^2(x)\,dx\Big]\overset{\mathcal{L}}{\longrightarrow} N(0,4\sigma^2),
\]
and if we further assume that the kernel function  $K$ is nonnegative, the estimator $-\ln(\frac{1}{n}+T_n(h_n))$ of the quadratic R\'{e}nyi entropy $R(f)=-\ln(\int_{\mathbb{R}}f^2(x)dx)$ satisfies 
\[
\sqrt{n}\Big[-\ln(\frac{1}{n}+T_n(h_n))-R(f)\Big]\overset{\mathcal{L}}{\longrightarrow} N\Big(0,\frac{4\sigma^2}{(\int_{\mathbb{R}}f^2(x)dx)^2}\Big)
\]
for some $\sigma^2\in (0,\infty)$.

\end{corollary}

Throughout this paper, if not mentioned otherwise, the letter $c$ with or without a subscript denotes a generic positive finite constant whose exact value is independent of $n$ and may change from line to line. For a complex number $z$, we use $\overline{z}$ and $|z|$ to denote its conjugate and modulus, respectively.  Moreover, we let $\phi(\lambda)$ be the characteristic function of linear process $X=\{X_n:\, n\in\mathbb{N}\}$. That is, $\phi(\lambda)=\mathbb{E}[e^{\iota \lambda X_n}]$ for each $\lambda\in\mathbb{R}$.

\section{Proof of the main result}

In this section, we will prove Theorem \ref{thm1}.  Firstly, we give two important definitions. For each $i\in\mathbb{Z}$, let $\mathcal{F}_i$ be the $\sigma$-field generated by $\{\varepsilon_k: k\leq i\}$. Given an integrable complex-valued random variable $Z$, we define the following projection operator $\mathscr{P}_i$ as
\[
\mathscr{P}_i Z=\mathbb{E}[Z|\mathcal{F}_i]-\mathbb{E}[Z|\mathcal{F}_{i-1}]
\]
for each $i\in\mathbb{Z}$. Then, for any two integrable complex-valued random variables $Z$ and $W$, $\mathbb{E}[\mathscr{P}_i Z\mathscr{P}_jW]=0$ if $i\neq j$.

For any $\delta\in(0,1)$, define
\begin{align} \label{tdeltan}
T^{\delta}_n(h_n)=\frac{2}{n(n-1)h_n} \sum_{1\le j< i\le n}K_{\delta}\Big(\frac{X_i-X_j}{h_n}\Big),
\end{align}
where $K_{\delta}(x)=\int_{\mathbb{R}}K(y)\frac{1}{\sqrt{2\pi}\delta}e^{-\frac{(x-y)^2}{2\delta^2}}dy$. 

It is easy to see that $K_{\delta}(x)$ is a symmetric and bounded function with $\int_{\R} K_{\delta}(u)\, du = 1$ and 
\begin{align} \label{2m}
&\int_{\R} u^2|K_{\delta}(u)|\, du \nonumber \\  \nonumber
&\leq \int_{\R} \int_{\R}u^2 |K(y)|\frac{1}{\sqrt{2\pi}\delta}e^{-\frac{(u-y)^2}{2\delta^2}}dy\, du \nonumber \\
&\leq 2\int_{\R} \int_{\R}(u-y)^2|K(y)|\frac{1}{\sqrt{2\pi}\delta}e^{-\frac{(u-y)^2}{2\delta^2}}du\, dy+2\int_{\R} \int_{\R}y^2|K(y)|\frac{1}{\sqrt{2\pi}\delta}e^{-\frac{(u-y)^2}{2\delta^2}}du\, dy \nonumber\\
&\leq 2\int_{\R}(1+|y|^2)|K(y)|dy<\infty.
\end{align}
Moreover, for any $1\leq i<j\leq n$, by the bounded convergence theorem,
\begin{align*}
\lim_{\delta\downarrow 0}\mathbb{E}\left|K\Big(\frac{X_i-X_j}{h_n}\Big)-K_{\delta}\Big(\frac{X_i-X_j}{h_n}\Big)\right|^2=0.
\end{align*}
Therefore, $\limsup\limits_{\delta\downarrow 0}\mathbb{E}|T_n(h_n)-T^{\delta}_n(h_n)|^2=0$.

\noindent
{\it Proof of Theorem \ref{thm1}} The proof will be done in several steps. 

\medskip
\noindent
{\bf Step 1.} We give the estimation for $\big|\E T_n(h_n)-\int_{\mathbb{R}}f^2(x)\mathrm{d}x\big|$. Using the Fourier inverse transform, Plancherel formula, symmetry of $K_{\delta}(x)$ and $\int_{\mathbb{R}}K_{\delta}(u)du=1$, 
\begin{align*}
\E T^{\delta}_n(h_n)-\int_{\R} f^2(x)\, dx
&=\frac{1}{\pi n(n-1)}\sum_{1\leq i<j\leq n}\int_{\R}\widehat{K_{\delta}}(\lambda h_n)\, \E\big[e^{-\iota \lambda (X_i-X_j)}\big]\, d\lambda-\frac{1}{2\pi}\int_{\R}|\phi(\lambda)|^2\, d\lambda\\
&=\frac{1}{\pi n(n-1)}\sum_{1\leq i<j\leq n}\int_{\R}\widehat{K_{\delta}}(\lambda h_n) \, \E\big[ (e^{-\iota \lambda X_i}-\phi(-\lambda))(e^{\iota \lambda  X_j}-\phi(\lambda))\big]\, d\lambda\\
&\qquad\qquad\qquad+\frac{1}{2\pi}\int_{\R}\big(\widehat{K_{\delta}}(\lambda h_n)-\widehat{K_{\delta}}(0)\big)|\phi(\lambda)|^2\, d\lambda\\
&=: \text{I}_1+\text{I}_2.
\end{align*}
By Lemma 7.1 in \cite{ssx}, $\int_{\mathbb{R}}|\lambda|^{2\gamma}|\phi_{\varepsilon}(\lambda)|^2<\infty$ and
\[
|\widehat{K_{\delta}}(\lambda)|=|\widehat{K}(\lambda)e^{-\frac{1}{2}\delta^2\lambda^2}|\leq |\widehat{K}(\lambda)|\leq \int_{\R}|K(x)|dx<\infty
\]
for all $\lambda\in\R$, we can obtain that $\text{I}_1$ is less than a constant multiple of $\frac{1}{n}$. 

Moreover, by symmetry of $K_{\delta}(x)$ and $\int_{\mathbb{R}}K_{\delta}(u)du=1$,
\begin{align}\label{kd0}
\Big|\widehat{K_{\delta}}(\lambda h_n)-\widehat{K}_{\delta}(0)\Big|=\Big|\int_{\R} (e^{\iota\lambda h_n u}-1)K_{\delta}(u) \, du\Big|\leq |\lambda|^{2\gamma} h^{2\gamma}_n \int_{\R} u^2 |K_{\delta}(u)|\, du.
\end{align} 
Then, by $\int_{\mathbb{R}}|\lambda|^{2\gamma}|\phi_{\varepsilon}(\lambda)|^2<\infty$ and (\ref{2m}), $|\text{I}_2|$ is less than a constant multiple of $h^{2\gamma}_n$. 

Hence,
\begin{align*}
\Big| \E T_n(h_n)-\int_{\R} f^2(x)\, dx \Big|
&\leq \limsup_{\delta\downarrow 0}\Big| \E T_n(h_n)-\E T^{\delta}_n(h_n) \Big|+\limsup_{\delta\downarrow 0}\Big| \E T^{\delta}_n(h_n)-\int_{\R} f^2(x)\, dx \Big|\\
&\leq c_1\Big(\frac{1}{n}+h^{2\gamma}_n\Big).
\end{align*}

\noindent
{\bf Step 2.} We give the decomposition for $T^{\delta}_n(h_n)-\E T^{\delta}_n(h_n)$. 

Recall the definition of $T^{\delta}_n(h_n)$ in (\ref{tdeltan}). Then, by the Fourier inverse transform,
\begin{align} \label{decomp}
T^{\delta}_n(h_n)-\E T^{\delta}_n(h_n)=A_{n,\delta}+D_{n,\delta}-\E D_{n,\delta}+N_{n,\delta},
\end{align}
where
\begin{align*}
A_{n,\delta}&=\frac{2}{n(n-1)h_n}\sum_{1\leq i<j\leq n,\, j-i\leq p_2} \Big[K_{\delta}\big(\frac{X_i-X_j}{h_n}\big)-\E K_{\delta}\big(\frac{X_i-X_j}{h_n}\big)\Big],\\
D_{n,\delta}&=\frac{1}{\pi n(n-1)}\sum_{1\leq i<j\leq n,\,  j-i>p_2}\int_{\R} \widehat{K_{\delta}}(\lambda h_n)\big(e^{-\iota \lambda  X_i}-\phi(-\lambda) \big)\big(e^{\iota \lambda  X_j}-\phi(\lambda) \big) \, d\lambda,\\
N_{n,\delta}&=\frac{1}{\pi n(n-1)}\sum_{1\leq i<j\leq n,\, j-i>p_2}\int_{\R} \widehat{K_{\delta}}(\lambda h_n) \Big[\big(e^{-\iota \lambda X_i}-\phi(-\lambda) \big) \phi(\lambda)+\big(e^{\iota \lambda X_j}-\phi(\lambda) \big) \phi(-\lambda)\Big] \, d\lambda
\end{align*}
and $p_2=\min\{i>p_1: a_i\neq 0\}$ with $p_1=\min\{i\geq 1: a_i\neq 0\}$. The finiteness of $p_1$ and $p_2$ follow from the assumption {\bf (A1)}.

\medskip
\noindent
{\bf Step 3.} We estimate $\E|A_{n,\delta}|^2$.  For any $1\leq i<j\leq n$, define
\[
A^{i,j}_{n,\delta}=\frac{1}{h_n}\Big[K_{\delta}\big(\frac{X_i-X_j}{h_n}\big)-\E K_{\delta}\big(\frac{X_i-X_j}{h_n}\big)\Big].
\]
Then
\begin{align} \label{and1}
\E|A_{n,\delta}|^2
&=\frac{4}{n^2(n-1)^2}\E \Big|\sum_{1\leq i<j\leq n,\, j-i\leq p_2} A^{i,j}_{n,\delta}\Big|^2 \nonumber\\
&=\frac{4}{n^2(n-1)^2}\sum_{\substack{1\leq i_{\theta}<j_{\theta}\leq n,\, j_{\theta}-i_{\theta}\leq p_2\\ \theta=1,2}} \E\big[A^{i_1,j_1}_{n,\delta} A^{i_2,j_2}_{n,\delta}\big]\nonumber\\
&\leq \frac{c_2}{n^3h^2_n}+\frac{c_2}{n^4} \sum_{\substack{1\leq i_{\theta}<j_{\theta}\leq n,\, j_{\theta}-i_{\theta}\leq p_2\\ \theta=1,2,\,|j_1-j_2|>4p_2}} \Big|\E\big[A^{i_1,j_1}_{n,\delta} A^{i_2,j_2}_{n,\delta}\big]\Big|\nonumber\\
&\leq \frac{c_2}{n^3h^2_n}+\frac{c_2}{n^4} \sum_{\substack{1\leq i_{\theta}<j_{\theta}\leq n,\, j_{\theta}-i_{\theta}\leq p_2\\ \theta=1,2,\, |j_1-j_2|>4p_2}} \sum^{\min\{j_1,j_2\}}_{k=-\infty}\Big|\E\big[\mathscr{P}_k(A^{i_1,j_1}_{n,\delta})\mathscr{P}_k(A^{i_2,j_2}_{n,\delta})\big]\Big|.
\end{align}

For any $1\leq i<j\leq n$, it is easy to see that $X_i-X_j=\sum\limits^{\infty}_{\ell=0} a^{i,j}_{\ell} \varepsilon_{j-\ell}$ where 
\begin{align} \label{aijl}
a^{i,j}_{\ell}=\left\{\begin{array}{ll}
-a_{\ell}  &  \text{if}\;  0\leq \ell<j-i,\\ \\
a_{\ell-j+i}-a_{\ell}  &    \text{if}\; \ell\geq j-i.
\end{array}
\right.
\end{align}

By the Fourier inverse transform,
\begin{align*}
\mathscr{P}_k(A^{i,j}_{n,\delta})
&=\mathbb{E}\Big[K_{\delta}\big(\frac{X_i-X_j}{h_n}\big)|\mathcal{F}_k\Big]-\mathbb{E}\Big[K_{\delta}\big(\frac{X_i-X_j}{h_n}\big)|\mathcal{F}_{k-1}\Big]\\
&=\frac{1}{2\pi}\int_{\mathbb{R}}\widehat{K_{\delta}}(\lambda h_n)\Big(\prod^{j-k-1}_{\ell=0} \phi_{\varepsilon}(-\lambda a^{i,j}_{\ell})\Big)\big(e^{-\iota\lambda a^{i,j}_{j-k}\varepsilon_{k}}-\phi_{\varepsilon}(-\lambda a^{i,j}_{j-k})\big)e^{-\iota \lambda \sum\limits^{\infty}_{\ell=j-k+1} a^{i,j}_{\ell}\varepsilon_{j-\ell}}\, d\lambda.
\end{align*}

For the case $k<j_1<j_2$, by Cauchy-Schwarz inequality and assumption {\bf (A2)}, 
\begin{align*}
&\Big|\E[\mathscr{P}_k(A^{i_1,j_1}_{n,\delta})\mathscr{P}_k(A^{i_2,j_2}_{n,\delta})]\Big|\\
&\leq \int_{\R^2}|\widehat{K}(\lambda_1 h_n)||\widehat{K}(\lambda_2 h_n)||\phi_{\varepsilon}(a_0\lambda_1)||\phi_{\varepsilon}(a_0\lambda_2)| \\
&\qquad\qquad \times \Big|\E\big[(e^{-\iota \lambda_1 a^{i_1,j_1}_{j_1-k} \varepsilon_k}-\phi_{\varepsilon}(-\lambda_1 a^{i_1,j_1}_{j_1-k}))(e^{-\iota \lambda_1 a^{i_2,j_2}_{j_2-k} \varepsilon_k}-\phi_{\varepsilon}(-\lambda_2 a^{i_2,j_2}_{j_2-k}))\big]\Big|  d\lambda_1\, d\lambda_2\\
&\leq c_3\, |a^{i_1,j_1}_{j_1-k}|^{\gamma}|a^{i_2,j_2}_{j_2-k}|^{\gamma} \int_{\R^2}|\widehat{K}(\lambda_1 h_n)||\widehat{K}(\lambda_2 h_n)||\phi_{\varepsilon}(a_0\lambda_1)||\phi_{\varepsilon}(a_0\lambda_2)||\lambda_1|^{\gamma}|\lambda_2|^{\gamma} d\lambda_1\, d\lambda_2\\
&=c_3\, |a^{i_1,j_1}_{j_1-k}|^{\gamma}|a^{i_2,j_2}_{j_2-k}|^{\gamma} \left(\int_{\R}|\widehat{K}(\lambda h_n)||\phi_{\varepsilon}(a_0\lambda)||\lambda|^{\gamma}\, d\lambda \right)^2\\
&\leq c_3\, |a^{i_1,j_1}_{j_1-k}|^{\gamma}|a^{i_2,j_2}_{j_2-k}|^{\gamma}\, \int_{\R}|\widehat{K}(\lambda h_n)|^2\, d\lambda\, \int_{\R}|\phi_{\varepsilon}(a_0\lambda)|^2|\lambda|^{2\gamma}\, d\lambda\\
&\leq c_4\, h^{-1}_n\,|a^{i_1,j_1}_{j_1-k}|^{\gamma}|a^{i_2,j_2}_{j_2-k}|^{\gamma}.
\end{align*}

Recall  the definition of $a^{i,j}_{\ell}$ in (\ref{aijl}) and $\sum\limits^{\infty}_{i=0} |a_i|^{\gamma}<\infty$. If $j-i=p_1$ and $a_0=a_{p_1}$, then there exists $\ell_0\geq p_1$ such that $a^{i,j}_{\ell_0}\neq 0$.  Now, for the case $k=j_1<j_2$,  by Cauchy-Schwarz inequality and assumption {\bf (A2)},
\begin{align*}
&\Big|\E[\mathscr{P}_{j_1}(A^{i_1,j_1}_{n,\delta})\mathscr{P}_{j_1}(A^{i_2,j_2}_{n,\delta})]\Big|\\
&\leq \int_{\R^2}|\widehat{K}(\lambda_1 h_n)||\widehat{K}(\lambda_2 h_n)||\phi_{\varepsilon}(a_0\lambda_2)||\phi_{\varepsilon}(a^{i_1,j_1}_p\lambda_1- a^{i_2,j_2}_{j_2-j_1+p}\lambda_2)| \\
&\qquad\qquad\times \bigg|\E\Big[\big(e^{-\iota \lambda_1 a^{i_1,j_1}_0 \varepsilon_{j_1}}-\phi_{\varepsilon}(-\lambda_1 a^{i_1,j_1}_0)\big)\big(e^{-\iota \lambda_2 a^{i_2,j_2}_{j_2-j_1} \varepsilon_{j_1}}-\phi_{\varepsilon}(-\lambda_2 a^{i_2,j_2}_{j_2-j_1})\big)\Big]\bigg|\,  d\lambda_1\, d\lambda_2\\
&\leq c_5\,|a^{i_2,j_2}_{j_2-j_1}|^{\gamma} \int_{\R}|\widehat{K}(\lambda_2 h_n)||\phi_{\varepsilon}(a_0\lambda_2)| |\lambda_2|^{\gamma}\Big(\int_{\R}|\widehat{K}(\lambda_1 h_n)||\phi_{\varepsilon}(a^{i_1,j_1}_p\lambda_1- a^{i_2,j_2}_{j_2-j_1+p}\lambda_2)|\, d\lambda_1\Big)\, d\lambda_2\\
 &\leq c_6\, h^{-1}_n\, |a^{i_2,j_2}_{j_2-j_1}|^{\gamma},
\end{align*}
where \[
a^{i_1,j_1}_p=\left\{
\begin{array}{cc}
-a_{p_1}   &   \text{if}\; p_1<j_1-i_1\\
a_0    &   \text{if}\; p_1>j_1-i_1\\
a_0-a_{p_1}     &   \text{if}\; j_1-i_1=p_1\; \text{and}\; a_0\neq a_{p_1}\\
a^{i_1,j_1}_{\ell_0}  &   \text{if}\; j_1-i_1=p_1\; \text{and}\; a_0=a_{p_1}
\end{array}
\right.
\] with
\[
p=\left\{
\begin{array}{cc}
j_1-p_1    &   \text{if}\; p_1<j_1-i_1\\
j_1-i_1    &   \text{if}\; p_1>j_1-i_1\\
p_1       &   \text{if}\; j_1-i_1=p_1\; \text{and}\; a_0\neq a_{p_1}\\
\ell_0  &   \text{if}\; j_1-i_1=p_1\; \text{and}\; a_0=a_{p_1}
\end{array}
\right..
\]
Therefore, by the definition of $a^{i,j}_{\ell}$ in (\ref{aijl}) and assumption {\bf (A1)},
\begin{align} \label{and2}
&\sum_{\substack{ 1\leq i_{\theta}<j_{\theta}\leq n, \, j_{\theta}-i_{\theta}\leq p_2\\
\theta=1,2,\,|j_1-j_2|>4p_2,\, j_1<j_2}}\sum^{\min\{j_1,j_2\}}_{k=-\infty}\Big|\E\big[\mathscr{P}_k(A^{i_1,j_1}_{n,\delta})\mathscr{P}_k(A^{i_2,j_2}_{n,\delta})\big]\Big| \nonumber\\
&\leq c_7\, h^{-1}_n \sum_{\substack{ 1\leq i_{\theta}<j_{\theta}\leq n, \, j_{\theta}-i_{\theta}\leq p_2\\
\theta=1,2,\, |j_1-j_2|>4p_2,\,j_1<j_2}}\Big(|a^{i_2,j_2}_{j_2-j_1}|^{\gamma}+\sum^{j_1-1}_{k=-\infty} |a^{i_1,j_1}_{j_1-k}|^{\gamma}|a^{i_2,j_2}_{j_2-k}|^{\gamma}\Big) \nonumber\\
&\leq c_8\, n h^{-1}_n.
\end{align}
Similarly,
\begin{align} \label{and3}
&\sum_{\substack{ 1\leq i_{\theta}<j_{\theta}\leq n, j_{\theta}-i_{\theta}\leq p_2\\
\theta=1,2,\, |j_1-j_2|>4p_2,\, j_1>j_2}}\sum^{\min\{j_1,j_2\}}_{k=-\infty}\Big|\E\big[\mathscr{P}_k(A^{i_1,j_1}_{n,\delta})\mathscr{P}_k(A^{i_2,j_2}_{n,\delta})\big]\Big|\leq c_9\, n h^{-1}_n.
\end{align}

Combining (\ref{and1}), (\ref{and2}) and (\ref{and3}) gives 
\begin{align} \label{ande}
\E|A_{n,\delta}|^2\leq \frac{c_{10}}{n^3 h^2_n}.
\end{align}

\medskip
{\bf Step 4.} We estimate $\E|D_{n,\delta}|$.  Note that 
\begin{align}\label{dnd}
D_{n,\delta}=D^1_{n,\delta}+D^2_{n,\delta},
\end{align}
where 
\begin{align*}
D^1_{n,\delta}=\frac{1}{\pi n(n-1)}\sum_{1\leq i<j\leq n,\,  j-i>p_2}\sum^i_{k=-\infty}\int_{\R} \widehat{K_{\delta}}(\lambda h_n)\mathscr{P}_k\Big(H(X_i)(-\lambda)\Big)\mathscr{P}_k\Big(H(X_j)(\lambda)\Big) \, d\lambda
\end{align*}
and 
\begin{align*}
D^2_{n,\delta}=\frac{1}{\pi n(n-1)}\sum_{\substack{1\leq i<j\leq n,\, j-i>p_2\\ k\leq i,\, \ell\leq j,\, k\neq \ell}}\int_{\R} \widehat{K_{\delta}}(\lambda h_n) \mathscr{P}_k\Big(H(X_i)(-\lambda)\Big)\mathscr{P}_{\ell}\Big(H(X_j)(\lambda)\Big) \, d\lambda
\end{align*}
with $H(X_i)(\lambda)=e^{\iota \lambda X_i}-\E e^{\iota \lambda X_i}$ for $1\leq i\leq n$ and $\lambda\in\R$.

Using similar arguments as in {\bf Step 3.},
\begin{align} \label{d1nd}
\E|D^1_{n,\delta}|
&=\frac{2}{\pi n(n-1)}\E\bigg|\sum_{1\leq i<j\leq n,\,  j-i>p_2}\sum^i_{k=-\infty}\int_{\R} \widehat{K_{\delta}}(\lambda h_n) \mathscr{P}_k\Big(H(X_i)(-\lambda)\Big)\mathscr{P}_k\Big(H(X_j)(\lambda)\Big) \, d\lambda\bigg| \nonumber\\
&\leq \frac{c_{11}}{n^2}\sum_{1\leq i<j\leq n,\,  j-i> p_2} \sum^{i}_{k=-\infty} |a_{i-k}|^{\gamma}|a_{j-k}|^{\gamma} \int_{\R}  |\phi_{\varepsilon}(a_0\lambda)||\phi_{\varepsilon}(a_{p_1}\lambda)||\lambda|^{2\gamma} \, d\lambda \nonumber \\
&\leq \frac{c_{12}}{n}.
\end{align}
Moreover,
\begin{align}\label{d2nd}
D^2_{n,\delta}=\text{II}_1+ \text{II}_2+\text{II}_3+\text{II}_4+\text{II}_5,
\end{align}
where
\begin{align*}
\text{II}_1&=\frac{1}{\pi n(n-1)}\sum_{1\leq i<j\leq n,\, j-i>p_2}\int_{\R} \widehat{K_{\delta}}(\lambda h_n) \mathscr{P}_i\Big(H(X_i)(-\lambda)\Big)\mathscr{P}_{j}\Big(H(X_j)(\lambda)\Big) \, d\lambda,\\
\text{II}_2&=\frac{1}{\pi n(n-1)}\sum_{1\leq i<j\leq n,\, j-i>p_2}\int_{\R} \widehat{K_{\delta}}(\lambda h_n) \mathscr{P}_i\Big(H(X_i)(-\lambda)\Big)\mathscr{P}_{j-p_1}\mathscr{P}_{j}\Big(H(X_j)(\lambda)\Big) \, d\lambda,\\
\text{II}_3&=\frac{1}{\pi n(n-1)}\sum_{1\leq i<j\leq n,\, j-i>p_2}\int_{\R} \widehat{K_{\delta}}(\lambda h_n) \mathscr{P}_{i-p_1}\Big(H(X_i)(-\lambda)\Big)\mathscr{P}_{j}\Big(H(X_j)(\lambda)\Big) \, d\lambda,\\
\text{II}_4&=\frac{1}{\pi n(n-1)}\sum_{1\leq i<j\leq n,\, j-i>p_2}\int_{\R} \widehat{K_{\delta}}(\lambda h_n) \mathscr{P}_{i-p_1}\Big(H(X_i)(-\lambda)\Big)\mathscr{P}_{j-p_1}\Big(H(X_j)(\lambda)\Big) \, d\lambda
\end{align*}
and
\begin{align*}
\text{II}_5&=\frac{1}{\pi n(n-1)}\sum_{\substack{1\leq i<j\leq n,\, j-i>p_2,\, k\neq \ell\\
i-k>p_1\, \text{or}\, j-\ell>p_1}}\int_{\R} \widehat{K_{\delta}}(\lambda h_n) \mathscr{P}_{k}\Big(H(X_i)(-\lambda)\Big)\mathscr{P}_{\ell}\Big(H(X_j)(\lambda)\Big)\, d\lambda.
\end{align*}

By Cauchy-Schwarz inequality and assumptions {\bf (A1)}, {\bf (A2)} and {\bf (A3)},
\begin{align}\label{II1}  
\mathbb{E}|\text{II}_{1}|^2 
&=\frac{1}{\pi^2 n^2(n-1)^2}\sum_{\substack{1\leq i_{\theta}<j\leq n, \, j-i_{\theta}>p_2\nonumber\\
 \theta=1,2}}\int_{\R^2} \widehat{K_{\delta}}(\lambda_1 h_n)\widehat{K_{\delta}}(\lambda_2 h_n) \nonumber\\
&\quad\times  \mathbb{E}\bigg[\mathscr{P}_{i_1}\big(H(X_{i_1})(-\lambda_1)\big) \mathscr{P}_{j}\big(H(X_{j})(\lambda_1)\big)  \mathscr{P}_{i_2}\big(H(X_{i_2})(-\lambda_2)\big)  \mathscr{P}_{j}\big(H(X_{j})(\lambda_2)\big) \bigg]\, d\lambda_1d\lambda_2 \nonumber\\
&\leq \frac{c_{13}}{n^4}\sum_{\substack{1\leq i_{\theta}<j\leq n, \, j-i_{\theta}>p_2\\ \theta=1,2}}\int_{\R^2} \Big|\widehat{K}(\lambda_1 h_n)\widehat{K}(\lambda_2 h_n) \phi_{\varepsilon}\big(a_{p_1}(\lambda_1+\lambda_2)\big)\phi_{\varepsilon}\big(a_{p_2}(\lambda_1+\lambda_2)\big)\Big| \nonumber\\
&\qquad\times\bigg\{ \bigg|\mathbb{E}\Big[\big(e^{\iota a_{j-i_1}(\lambda_1+\lambda_2)}-\phi_{\varepsilon}(a_{j-i_1}(\lambda_1+\lambda_2))\big)\big(e^{-\iota a_0\lambda_1}-\phi_{\varepsilon}(-a_0\lambda_1)\big)\Big]\bigg|1_{\{i_1>i_2\}}\nonumber\\
&\qquad\qquad\quad+\bigg|\mathbb{E}\Big[\big(e^{\iota a_{j-i_2}(\lambda_1+\lambda_2)}-\phi_{\varepsilon}(a_{j-i_2}(\lambda_1+\lambda_2))\big)(e^{-\iota a_0\lambda_1}-\phi_{\varepsilon}(-a_0\lambda_1)\Big]\bigg|1_{\{i_2>i_1\}}\nonumber\\
&\qquad\qquad\qquad\quad+1_{\{i_2=i_1\}}\bigg\}\, d\lambda_1d\lambda_2\nonumber\\
&\leq \frac{c_{14}}{n^4}\sum_{\substack{1\leq i_{\theta}<j\leq n, \, j-i_{\theta}>p_2\nonumber\\
\theta=1,2}}\int_{\R^2} \Big|\widehat{K}(\lambda_1 h_n)\widehat{K}(\lambda_2 h_n) \phi_{\varepsilon}\big(a_{p_1}(\lambda_1+\lambda_2)\big)\phi_{\varepsilon}\big(a_{p_2}(\lambda_1+\lambda_2)\big)\Big| \nonumber\\
&\qquad\times\Big\{|\lambda_1+\lambda_2|^{\gamma}\Big(|a_{j-i_1}|^{\gamma}1_{\{i_1>i_2\}}+|a_{j-i_2}|^{\gamma}1_{\{i_2>i_1\}}\Big)+1_{\{i_2=i_1\}}\Big\}\, d\lambda_1\, d\lambda_2\nonumber\\
&\leq \frac{c_{15}}{n^2h_n}
\end{align}
and
\begin{align} \label{II2}  
\mathbb{E}|\text{II}_{2}|^2
&=\frac{1}{\pi^2 n^2(n-1)^2}\sum_{\substack{1\leq i_{\theta}<j\leq n, \, j-i_{\theta}>p_2 \nonumber\\
 \theta=1,2}}\int_{\R^2} \widehat{K_{\delta}}(\lambda_1 h_n)\widehat{K_{\delta}}(\lambda_2 h_n) \nonumber\\
&\quad\times  \mathbb{E}\bigg[\mathscr{P}_{i_1}\big(H(X_{i_1})(-\lambda_1)\big) \mathscr{P}_{j-p_1}\big(H(X_{j})(\lambda_1)\big)  \mathscr{P}_{i_2}\big(H(X_{i_2})(-\lambda_2) \big) \mathscr{P}_{j-p_1}\big(H(X_{j})(\lambda_2)\big) \bigg]\, d\lambda_1\,d\lambda_2\nonumber\\
&\leq \frac{c_{16}}{n^4}\sum_{\substack{1\leq i_{\theta}<j\leq n, \, j-i_{\theta}>p_2\\ \theta=1,2}}\int_{\R^2} \Big|\widehat{K}(\lambda_1 h_n)\widehat{K}(\lambda_2 h_n) \phi_{\varepsilon}(a_0\lambda_1)\phi_{\varepsilon}(a_0\lambda_2)\phi_{\varepsilon}\big(a_{p_2}(\lambda_1+\lambda_2)\big)\Big| \nonumber\\
&\qquad\times\bigg\{ \bigg|\mathbb{E}\Big[\big(e^{\iota a_{j-i_1}(\lambda_1+\lambda_2)}-\phi_{\varepsilon}(a_{j-i_1}(\lambda_1+\lambda_2))\big)\big(e^{-\iota a_0\lambda_1}-\phi_{\varepsilon}(-a_0\lambda_1)\big)\Big]\bigg|1_{\{i_1>i_2\}}\nonumber\\
&\qquad\qquad\quad+\bigg|\mathbb{E}\Big[\big(e^{\iota a_{j-i_2}(\lambda_1+\lambda_2)}-\phi_{\varepsilon}(a_{j-i_2}(\lambda_1+\lambda_2))\big)(e^{-\iota a_0\lambda_1}-\phi_{\varepsilon}(-a_0\lambda_1)\Big]\bigg|1_{\{i_2>i_1\}}\nonumber\\
&\qquad\qquad\qquad\quad+1_{\{i_2=i_1\}}\bigg\}\, d\lambda_1\, d\lambda_2\nonumber\\
&\leq \frac{c_{17}}{n^4}\sum_{\substack{1\leq i_{\theta}<j\leq n, \, j-i_{\theta}>p_2\\ \theta=1,2}}\int_{\R^2} \Big|\widehat{K}(\lambda_1 h_n)\widehat{K}(\lambda_2 h_n) \phi_{\varepsilon}(a_0\lambda_1)\phi_{\varepsilon}(a_0\lambda_2)\phi_{\varepsilon}\big(a_{p_2}(\lambda_1+\lambda_2)\big)\Big| \nonumber\\
&\qquad\times\bigg\{|\lambda_1+\lambda_2|^{\gamma}\Big(|a_{j-i_1}|^{\gamma}1_{\{i_1>i_2\}}+|a_{j-i_2}|^{\gamma}1_{\{i_2>i_1\}}\Big)+1_{\{i_2=i_1\}}\bigg\}\, d\lambda_1\, d\lambda_2\nonumber\\
&\leq \frac{c_{18}}{n^2h_n}.
\end{align}
Similarly,
\begin{align} \label{II3}
\mathbb{E}|\text{II}_{3}|^2\leq \frac{c_{19}}{n^2h_n}\quad \text{and}\quad \mathbb{E}|\text{II}_{4}|^2\leq \frac{c_{19}}{n^2h_n}.
\end{align}
Note that 
\begin{align*}
\text{II}_4=\text{II}_{41}+\text{II}_{42},
\end{align*}
where 
\begin{align*}
\text{II}_{41}&=\frac{1}{\pi n(n-1)}\sum_{\substack{1\leq i<j\leq n,\,j-i>p_2,\, k>\ell\\
i-k>p_1\, \text{or}\, j-\ell>p_1}}\int_{\R} \widehat{K_{\delta}}(\lambda h_n) \mathscr{P}_i\Big(H(X_i)(-\lambda)\Big)\mathscr{P}_{\ell}\Big(H(X_j)(\lambda)\Big) \, d\lambda,
\end{align*}
and
\begin{align*}
\text{II}_{42}&=\frac{1}{\pi n(n-1)}\sum_{\substack{1\leq i<j\leq n,\,j-i>p_2, \, k<\ell\\
i-k>p_1\, \text{or}\, j-\ell>p_1}}\int_{\R} \widehat{K_{\delta}}(\lambda h_n) \mathscr{P}_{i}\Big(H(X_i)(-\lambda)\Big)\mathscr{P}_{\ell}\Big(H(X_j)(\lambda)\Big) \, d\lambda.
\end{align*}

Clearly,
\begin{align*}  
\mathbb{E}|\text{II}_{41}|^2
&=\frac{1}{\pi^2 n^2(n-1)^2}\sum_{\substack{1\leq i_{\theta}<j_{\theta}\leq n, \, j_{\theta}-i_{\theta}>p_2\\ i_{\theta}-k_{\theta}>p_1\,\text{or}\, j_{\theta}-\ell_{\theta}>p_1\\
k_{\theta}>\ell_{\theta},\,\theta=1,2}}\int_{\R^2} \widehat{K_{\delta}}(\lambda_1 h_n)\widehat{K_{\delta}}(\lambda_2 h_n) \, \text{I}^{i_1,i_2,j_1,j_2}_{k_1,k_2,\ell_1,\ell_2}(\lambda_1,\lambda_2)\, d\lambda_1d\lambda_2,
\end{align*}
where 
\begin{align} \label{expectation}
&\text{I}^{i_1,i_2,j_1,j_2}_{k_1,k_2,\ell_1,\ell_2}(\lambda_1,\lambda_2) \nonumber\\
&\qquad=\mathbb{E}\bigg[\mathscr{P}_{k_1}\Big(H(X_{i_1})(-\lambda_1)\Big) \mathscr{P}_{\ell_1}\Big(H(X_{j_1})(\lambda_1)\Big)  \mathscr{P}_{k_2}\Big(H(X_{i_2})(-\lambda_2)\Big)  \mathscr{P}_{\ell_2}\Big(H(X_{j_2})(\lambda_2)\Big) \bigg]
\end{align}
There are four possibilities for the orderings of $k_1, \ell_1, k_2,\ell_2$: 
\[
(1)\; k_1\geq k_2\geq \ell_1\geq \ell_2,  \; (2)\; k_1\geq k_2\geq \ell_2\geq \ell_1,\; (3)\; k_2\geq k_1\geq \ell_1\geq \ell_2, \; (4)\;  k_2\geq k_1\geq \ell_2\geq \ell_1.
\]
Note that the expectation in $\eref{expectation}$ is equal to zero if $k_1\neq k_2$. Hence, by Cauchy-Schwarz inequality and assumption {\bf (A2)},
\begin{align*}
&\Big|\text{I}^{i_1,i_2,j_1,j_2}_{k_1,k_2,\ell_1,\ell_2}(\lambda_1,\lambda_2) \Big|\\
&\leq 1_{\{k_1=k_2\}}\Big|\phi_{\varepsilon}(a_0\lambda_1)\phi_{\varepsilon}(a_{p_1}\lambda_1)\phi_{\varepsilon}(a_0\lambda_2)\phi_{\varepsilon}(a_{p_1}\lambda_2)\Big|\\
&\qquad\times \bigg|\mathbb{E}\Big[\big(e^{\iota a_{i_1-k_1}\lambda_1\varepsilon_{k_1}}-\phi_{\varepsilon}(a_{i_1-k_1}\lambda_1)\big)\big(e^{\iota a_{i_2-k_2}\lambda_2\varepsilon_{k_2}}-\phi_{\varepsilon}(a_{i_2-k_2}\lambda_2)\big)\Big]\bigg|\\
&\times \bigg\{\bigg|\mathbb{E}\Big[(e^{-\iota (a_{i_1-\ell_1}\lambda_1+a_{i_2-\ell_1}\lambda_2)\varepsilon_{\ell_1}}-\phi_{\varepsilon}(-a_{i_1-\ell_1}\lambda_1+a_{i_2-\ell_1}\lambda_2))(e^{\iota a_{j_1-\ell_1}\lambda_1\varepsilon_{\ell_1}}-\phi_{\varepsilon}(a_{j_1-\ell_1}\lambda_1))\Big]\bigg|1_{\{\ell_1>\ell_2\}}\\
&\qquad+\bigg|\mathbb{E}\Big[\big(e^{-\iota (a_{i_1-\ell_2}\lambda_1+a_{i_2-\ell_2}\lambda_2)\varepsilon_{\ell_2}}-\phi_{\varepsilon}(-a_{i_1-\ell_2}\lambda_1+a_{i_2-\ell_2}\lambda_2)\big)\big(e^{\iota a_{j_2-\ell_2}\lambda_2\varepsilon_{\ell_2}}-\phi_{\varepsilon}(a_{j_2-\ell_2}\lambda_2)\big)\Big]\bigg|1_{\{\ell_2>\ell_1\}}\\
&\qquad\qquad+\bigg|\mathbb{E}\Big[\big(e^{\iota a_{j_1-\ell_1}\lambda_1\varepsilon_{\ell_1}}-\phi_{\varepsilon}(a_{j_1-\ell_1}\lambda_1)\big)\big(e^{\iota a_{j_2-\ell_2}\lambda_2\varepsilon_{\ell_2}}-\phi_{\varepsilon}(a_{j_2-\ell_2}\lambda_2)\big)\Big]\bigg|1_{\{\ell_2=\ell_1\}}\bigg\}\\
&\leq c_{20}1_{\{k_1=k_2\}}\Big|\phi_{\varepsilon}(a_0\lambda_1)\phi_{\varepsilon}(a_{p_1}\lambda_1)\phi_{\varepsilon}(a_0\lambda_2)\phi_{\varepsilon}(a_{p_1}\lambda_2)\Big||a_{i_1-k_1}\lambda_1|^{\gamma}|a_{i_2-k_2}\lambda_2|^{\gamma}\\
&\qquad\times \Big\{|a_{i_1-\ell_1}\lambda_1+a_{i_2-\ell_1}\lambda_2|^{\gamma}|a_{j_1-\ell_1}\lambda_1|^{\gamma}||a_{j_2-\ell_2}\lambda_2|^{\gamma}|1_{\{\ell_1>\ell_2\}}\\
&\qquad\qquad\qquad+|a_{i_1-\ell_2}\lambda_1+a_{i_2-\ell_2}\lambda_2|^{\gamma}|a_{j_2-\ell_2}\lambda_2|^{\gamma}|a_{j_1-\ell_1}\lambda_1|^{\gamma}|1_{\{\ell_2>\ell_1\}}\\
&\qquad\qquad\qquad\qquad+|a_{j_1-\ell_1}\lambda_1|^{\gamma}|a_{j_2-\ell_2}\lambda_2|^{\gamma}1_{\{\ell_2=\ell_1\}}\Big\}.
\end{align*}  
Thus, by assumption (A3),
\begin{align*}  
\mathbb{E}|\text{II}_{41}|^2
&\leq \frac{c_{21}}{n^4}\sum_{\substack{1\leq i_{\theta}<j_{\theta}\leq n, \, j_{\theta}-i_{\theta}>p_2\\ i_{\theta}-k_{\theta}>p_1\,\text{or}\, j_{\theta}-\ell_{\theta}>p_1\\
k_{\theta}>\ell_{\theta},\, \theta=1,2}}\int_{\R^2} \Big|\widehat{K}(\lambda_1 h_n)\widehat{K}(\lambda_2 h_n)\Big|\Big|\text{I}^{i_1,i_2,j_1,j_2}_{k_1,k_2,\ell_1,\ell_2}(\lambda_1,\lambda_2)\Big|\, d\lambda_1\, d\lambda_2\leq \frac{c_{22}}{n^2}.
\end{align*}
Similarly, 
\[
\mathbb{E}|\text{II}_{42}|^2\leq \frac{c_{23}}{n^2}.
\] 
Therefore,
\begin{align} \label{II4}
\E|\text{II}_4|^2=\E|\text{II}_{41}+\text{II}_{42}|^2\leq \frac{c_{24}}{n^2}.
\end{align}

Combining (\ref{dnd}), (\ref{d1nd}), (\ref{d2nd}), (\ref{II1}), (\ref{II2}), (\ref{II3}) and (\ref{II4}) gives
\begin{align}  \label{dnde}
\E|D_{n,\delta}|\leq \frac{c_{25}}{\sqrt{n^2h_n}}.
\end{align}

\noindent
{\bf Step 5.} We estimate $\E\big|N_{n,\delta}-\overline{N}_n\big|^2$ where
\[
\overline{N}_n=\frac{1}{\pi n}\sum^n_{i=1}\int_{\R} \widehat{K}(0) \big(e^{-\iota \lambda X_i}-\phi(-\lambda) \big) \phi(\lambda) \,  d\lambda.
\]

Recall that
\begin{align*}
N_{n,\delta}&=\frac{1}{\pi n(n-1)}\sum_{1\leq i<j\leq n,\, j-i>p_2}\int_{\R} \widehat{K_{\delta}}(\lambda h_n) \Big[\big(e^{-\iota \lambda X_i}-\phi(-\lambda) \big) \phi(\lambda)+\big(e^{\iota \lambda X_j}-\phi(\lambda) \big) \phi(-\lambda)\Big] \, d\lambda\\
&=\frac{1}{2\pi n(n-1)}\sum_{1\leq i\neq j\leq n,\, |j-i|>p_2}\int_{\R} \widehat{K_{\delta}}(\lambda h_n) \Big[\big(e^{-\iota \lambda X_i}-\phi(-\lambda) \big) \phi(\lambda)+\big(e^{\iota \lambda X_j}-\phi(\lambda) \big) \phi(-\lambda)\Big] \, d\lambda,
\end{align*}
where we use the symmetry of $K_{\delta}(x)$ in the last equality. 

Let 
\[
\widetilde{N}_{n,\delta}=\frac{1}{\pi n}\sum^n_{i=1}\int_{\R} \widehat{K_{\delta}}(\lambda h_n) \big(e^{-\iota \lambda X_i}-\phi(-\lambda) \big) \phi(\lambda) \, d\lambda.
\] 
It is easy to see that $|N_{n,\delta}-\widetilde{N}_{n,\delta}|$ is less than a constant multiple of $\frac{1}{n}$. 

Moreover, by Cauchy-Schwarz inequality, Fubini theorem, (\ref{kd0}), assumption {\bf (A3)} and Lemma 7.2 in \cite{ssx},
\begin{align*}
\E\big|\widetilde{N}_{n,\delta}-\overline{N}_n\big|^2
&=\frac{1}{\pi^2}\E\left[ \Big|\int_{\R} \big(\widehat{K_{\delta}}(\lambda h_n)-\widehat{K_{\delta}}(0)\big)\Big(\frac{1}{n}\sum^n_{i=1} e^{-\iota \lambda X_i}-\phi(-\lambda) \Big) \phi(\lambda) \,  d\lambda\Big|^2\right]\\
&\leq \left(\int_{\R} \big|\widehat{K_{\delta}}(\lambda h_n)-\widehat{K_{\delta}}(0)\big|^2|\phi(\lambda)|\, d\lambda \right) \Big(\int_{\mathbb{R}} \mathbb{E} \big|\frac{1}{n}\sum^n_{i=1} e^{-\iota \lambda X_i}-\phi(\lambda)\big|^2\big|\phi(\lambda)\big| \,  d\lambda \Big)\\
&\leq c_{26} \frac{h^{2\gamma}_n}{n}.
\end{align*}
Hence 
\begin{align}  \label{Nne}
\E\big|N_n-\overline{N}_n\big|^2\leq c_{27} \Big(\frac{1}{n^2}+\frac{h^{2\gamma}_n}{n}\Big).
\end{align}
Note that $\widehat{K}(0)=1$. By the Fourier inverse transform and Plancherel formula,
\[
\overline{N}_n=\frac{2}{n} \sum\limits^n_{i=1}\Big(f(X_i)-\int_{\R} f^2(x)\, dx\Big).
\]  

Combining $\eref{decomp}$, $\eref{ande}$, $\eref{dnde}$ and $\eref{Nne}$ gives
\begin{align*}
&\E\Big|T_n(h_n)-\E T_n(h_n)-\frac{1}{n}\sum^n_{i=1}Y_i\Big|\\
&\leq \limsup_{\delta\downarrow 0}\E\left|T_n(h_n)-\E T_n(h_n)-\big(T^{\delta}_n(h_n)-\E T^{\delta}_n(h_n)\big)\right|+\limsup_{\delta\downarrow 0}\E\Big|T^{\delta}_n(h_n)-\E T^{\delta}_n(h_n)-\frac{1}{n}\sum^n_{i=1}Y_i\Big|\\
&\leq c_{28}\Big(\frac{1}{\sqrt{n^3h^2_n}}+\frac{1}{\sqrt{n^2h_n}}+\frac{h^{\gamma}_n}{\sqrt{n}}\Big).
\end{align*}

\medskip
\noindent
{\bf Step 6.} We show the central limit theorem for 
\[
\frac{2}{n} \sum\limits^n_{i=1}\Big(f(X_i)-\int_{\R} f^2(x)\, dx\Big).
\]
Using Lemma 1 in \cite{wu} and similar arguments as in the proof of Theorem 2.1 in \cite{ssx},
\[
\frac{1}{n} \sum\limits^n_{i=1}\Big(f(X_i)-\int_{\R} f^2(x)\, dx\Big)\overset{\mathcal{L}}{\longrightarrow} N(0,\sigma^2)
\]
for some $\sigma^2\in (0,\infty)$.

This finishes the proof of Theorem \ref{thm1}.
\qed

\bigskip

$\begin{array}{cc}
\begin{minipage}[t]{1\textwidth}
{\bf  Yudan Xiong}\\
School of Statistics, East China Normal University, Shanghai 200262, China \\
\texttt{xyd\_980107@163.com}
\end{minipage}
\hfill
\end{array}$

\medskip

$\begin{array}{cc}
\begin{minipage}[t]{1\textwidth}
{\bf Fangjun Xu}\\
KLATASDS-MOE, School of Statistics, East China Normal University, Shanghai, 200062, China \\
NYU-ECNU Institute of Mathematical Sciences at NYU Shanghai,  200062, China\\
\texttt{fjxu@finance.ecnu.edu.cn}
\end{minipage}
\hfill
\end{array}$

\end{document}